\documentclass{amsart}
\usepackage{amscd,amssymb,subfigure,hyperref, epsfig}
\usepackage[arrow,matrix,graph,frame,poly,arc,tips]{xy}
\usepackage{amsmath}

\theoremstyle{plain}
\newtheorem{thm}[subsection]{Theorem}

\newtheorem{lem}[subsection]{Lemma}

\newtheorem{cor}[subsection]{Corollary}

\theoremstyle{definition}

\newtheorem{example}[subsection]{Example}

\numberwithin{equation}{section}

\newcommand{\A}{{\mathcal A}}

\newcommand{\La}{L_{\A}}

\renewcommand{\k}{\Bbbk}
\newcommand{\TM}{D({\mathcal A})}
\newcommand{\TMp}{D({\mathcal A}_X)}
\newcommand{\TMG}{D({\mathcal A}_G)}

\newcommand{\tm}{D_0({\mathcal A})}

\newcommand{\rk}{\mathrm{rank}}
\newcommand{\pdim}{\mathrm{pdim}}

\begin{document}

\title{
Derivation modules of orthogonal duals of hyperplane arrangements}

\author[Joseph P.S. Kung]{Joseph P.S. Kung}
\address{Department of Mathematics,
University of North Texas,  Denton, TX 76203}
\email{\href{mailto: kung@unt.edu}{kung@unt.edu}}

\author[Hal Schenck]{Hal Schenck$^1$}
\address{Department of Mathematics,
Texas A\&M University, College Station, TX 77843}
\email{\href{mailto: schenck@math.tamu.edu}{schenck@math.tamu.edu}}
\urladdr{\href{http://www.math.tamu.edu/~schenck/}%
{http://www.math.tamu.edu/\~{}schenck}}

\thanks{$^1$Partially supported by NSF DMS 03-11142, NSA MDA
  904-03-1-0006, and ATP 010366-0103.}



\begin{abstract}
Let $A$ be an $n \times d$ matrix having full rank $n.$
An orthogonal dual $A^{\perp}$
of $A$ is a $(d-n) \times d$ matrix of rank $(d-n)$
such that every row of $A^{\perp}$ is orthogonal (under
the usual dot product) to every row of $A.$
We define the orthogonal dual for arrangements by identifying
an essential (central) arrangement of $d$ hyperplanes
in $n$-dimensional space with the $n \times d$ matrix of coefficients of
the homogeneous linear forms for which the hyperplanes are kernels.
When $n \ge 5,$ we show that if
the matroid (or the lattice of intersection) of
an $n$-dimensional essential arrangement $\A$
contains a modular copoint whose complement spans, then
the  derivation module of the
orthogonally dual arrangement $\A^{\perp}$
has projective dimension at
least $\lceil n(n+2)/4 \rceil - 3.$
\end{abstract}

\maketitle


\section{Introduction}
\label{sec:intro}
An important conjecture in the theory of hyperplane arrangements
is {\em Terao's conjecture} \cite{t}:
whether the derivation module $D(\A)$ of a central
arrangement $\A$ is free depends only on the ``combinatorics'',
that is to say, the matroid of $\A .$   Since being free is equivalent to
having zero projective dimension, a natural extension
of Terao's conjecture is that the projective dimension
$\pdim (D(\A))$ of
$D(\A)$ depends only on the matroid of $\A.$  No counterexamples
to this extended conjecture are known, although in \cite{z} 
Ziegler gives two arrangements with the same matroid but
 non-isomorphic derivation modules.  The
 common matroid in Ziegler's example has rank $3$ and is the truncation
 of the orthogonal dual $M^{\perp}(K_{3,3})$ of the
 cycle matroid of the complete bipartite graph $K_{3,3}.$
 This matroid has two inequivalent representations, one from projecting
 a representation of $M^{\perp}(K_{3,3})$ from a point
 in general position, the other,
 a less special one which cannot be ``erected''.  However,
 the derivation modules of these two arrangements both have projective
 dimension $1.$   Ziegler's examples provide one motivation
 to study derivation modules of orthogonal duals of arrangements.
 Another motivation comes from computer experiments suggesting that
 derivation modules of duals of free arrangements tend to
 have high projective dimension.

The notion of orthogonal duality is pervasive in combinatorics.
H. Whitney first defined duality for matroids
in his 1932 paper \cite{Whitney} to extend the notion of a
dual or face graph of a planar graph to arbitrary graphs.
He proved the theorem (equivalent to Kuratowski's theorem for planarity)
that a graph is planar if and only if its matroid dual is a graphic matroid.
Another example (suggested by a referee) is the concept of {\em association} introduced by Coble
in \cite{c}; for an application to generic arrangements see \cite{dk}.
Duality also occurs in linear programming, combinatorial optimization, 
and coding theory. It is closely  related to Alexandrov and other 
kinds of duality in algebraic topology. See, for example, \cite{CR,Crapo, OX}.

Let $X$ be a subspace in the lattice $\La$ of intersection of
the arrangement $\A.$  The {\em closed subarrangement}
$\A_X$ is the subset of all hyperplanes
in $\A$ containing $X.$  When $X$ is $1$-dimensional,
$\A_X$ is a {\em copoint}.
A closed subarrangement $\A_X$ (or its associated subspace $X$)
is {\em modular} if
$$
\rk(X \vee Y) + \rk(X \cap Y) = \rk(X) + \rk(Y)
$$
for every subspace $Y$ in $\La.$
Chains of modular flats occur (by definition)
in supersolvable arrangements.   In addition,
it is easy to show by induction and
the addition-deletion lemma (see \cite{t} or
\cite{ot}, Chap. 4) that if an arrangement $\A$
has a modular copoint $\A_X$ which is free, then $\A$ itself is free.

A subarrangement of an essential arrangement {\em spans} if it is essential. 
The {\em ceiling} $\lceil x \rceil$ of a 
real number $x$ is the the smallest integer greater than or equal to $x.$ The {\em floor}  
$\lfloor x \rfloor$ of $x$ is the largest integer less than or equal to $x.$ 

Our main result is the following theorem.

\begin{thm}
\label{thm:dualPD}
Let $\A$ be an essential arrangement over an arbitrary field
 with
 a modular copoint $X$ such that
 its complement $\A \backslash \A_X$ spans.
 Suppose that the dimension $n$ of $\A$ is at least $5.$  Then the
 projective dimension of the derivation module of the orthogonal dual
 $\A^{\perp}$ is bounded below by $\lceil n(n+2)/4 \rceil - 3.$
\end{thm}

\noindent 

The proof of Theorem \ref{thm:dualPD} is combinatorial.  We
show that an arrangement $\A$ satisfying the main hypotheses in
\ref{thm:dualPD} contains a spanning
subarrangement with the same matroid as the braid arrangement
$A_{n+1}.$
This implies that the dual $\A^{\perp}$ contains a closed circuit
with at least $\lceil n(n+2)/4 \rceil$ hyperplanes.  The proof 
is completed by combining a result of Terao on projective
dimension of closed subarrangements with results of
Rose and Terao, and Yuzvinsky
on the projective dimension of generic arrangements.

\section {Projective dimension of $\TM$ and closed subarrangements}

In this section, we discuss the two theorems from hyperplane arrangements
we need.  Both theorems hold over arbitrary fields.  

A {\em generic arrangement} is an arrangement of
at least $n+1$ hyperplanes in $n$-dimensional space for which
every subset of $n$ hyperplanes
is independent.  In particular, matroids of generic
arrangements are uniform matroids.
The following theorem is due to
Rose and Terao \cite{rt} and Yuzvinsky \cite{y}.

\begin{thm}
\label{thm:generic}
If $\A$ is a generic arrangement in $\k^n$, then $\pdim(\TM)=n-2.$
\end{thm}

We shall also use the following theorem of Terao \cite{t2} (see also \cite{bt}).

\begin{thm}
\label{thm:PD}
If $\A_X$ is a closed subarrangement of $\A,$ then
$$
\pdim(\TM) \ge \pdim(\TMp).
$$
\end{thm}

Terao's proof is unpublished.  Yuzvinsky gives a proof in \cite{y1}. For the reader's
convenience, we give another proof, which is a more elementary version
of the proof in \cite{y1} (but requires the hypothesis that the field $\k$ has characteristic zero).
Let $S$ be the symmetric algebra $\mathrm{Sym} (V^*)$ of the dual space $V^*.$
The algebra $S$ is isomorphic to the polynomial algebra
$\k [x_1,x_2, \ldots, x_n],$ where $\{x_i\}$ is a dual basis for $V.$
Let $\A = \{H_i: 1 \le i \le d\}$ and $Q$ be the polynomial
$ Q = \prod_{i=1}^d l_i, $
where for each $i,$
$l_i$ is a homogeneous linear form such that
the kernel $V(l_i)$ of $l_i$ is the hyperplane $H_i.$
The {\em derivation module} $\TM $ is the 
$S$-module of all $S$-derivations $\theta$ such that for all $i$, 
$\theta(l_i)$ is in the principal ideal $\langle l_i \rangle
\subseteq S$. If char $\k = 0$, this is equivalent to
the single condition $\theta(Q) \in \langle Q \rangle$.
The Euler derivation
$\sum  x_i \partial/\partial x_i$
generates a free summand $S(-1)$ of $\TM$ and
$$
\TM = S(-1) \oplus \tm,
$$
 where $\tm$ is the kernel of the Jacobian matrix $J_Q,$  the
 $n \times 1$ matrix with $i,1$-entry equal to
 $\partial Q/\partial x_i$ (see, for example, \cite{y}).
 In particular, the projective dimension of $\TM$ is one
 less than the projective dimension of the ideal
 $\langle  J_Q \rangle$ generated by
 the entries of the matrix $J_Q$.

Let $X$ be a subspace in the intersection lattice of $\A.$
Order the hyperplanes of $\A$ so that the closed subarrangement $\A_X$
equals $\{H_1, H_2, \ldots, H_s\}.$
Choose coordinates so that $X$ is the subspace $V(x_1,x_2,\ldots,x_k)$
defined by the equations $x_1 = 0, x_2 =0, \ldots, x_k = 0,$
and hence, a hyperplane $H_i$ in $\A_X$ may be written as
the kernel $V(l_i)$ with $l_i$ a homogeneous linear form in $\k[x_1,\ldots, x_k].$

Let $P$ be the prime 
ideal $\langle x_1,\ldots, x_k \rangle$ in $S.$  By our choice
of coordinates, if the hyperplane $V(l_i)$
does not contain the subspace $X,$ then $l_i$ equals
$\gamma_i + \delta_i$ where $\gamma_i$ is a linear form in
$\k[x_1,\ldots, x_k],$  $\delta_i$ is a linear form in
$\k[x_{k+1},\ldots, x_n],$ and $\delta_i \not\equiv 0.$

Write $Q = L  K$, where
$L = \prod_{i=1}^s l_i$ and $K = \prod_{i=s+1}^{|\A|} l_i.$
Computing the $i,1$-entry of $J_Q$ by the product rule, we have
$$
\frac {\partial Q}{\partial x_i} =
L \, \frac {\partial K}{\partial x_i} \,+\, K \, \frac {\partial L}{\partial x_i}.
$$
\noindent
By our choice of coordinates, $\partial L / \partial x_i  =0$ when
$i > k.$   Hence, the Jacobian matrix $J_Q$ simplifies to the transpose
of the matrix
$$
\left[L \, \frac{\partial K}{\partial x_1}  +  K \, \frac{\partial L}{\partial x_1},
 \ldots, L \, \frac{\partial K}{\partial x_k}
 +  K \, \frac{\partial L}{\partial x_k},
 L \, \frac{\partial K}{\partial x_{k+1}},
 \ldots, L \,\frac{\partial K}{\partial x_n} \right].
$$
We localize at the prime ideal $P.$
In the local ring $S_P,$ every element not in $P$ is a unit.
Since each  $l_i$ with  $i > s$ has the non-zero form
$\delta_i$ in $\k[x_{k+1},\ldots, x_n],$ the product $K$
contains at least one monomial in $\k[x_{k+1},\ldots, x_n].$
Hence $K$
is a unit in $S_P.$   Similarly, $\partial  K / \partial x_i$ is
nonzero for some $i \in \{s+1,\ldots,n \}$
{\em and} still contain a nonzero monomial in $\k[x_{k+1},\ldots, x_n].$
In particular, as an element in $S_P,$ $L$ equals $K^{-1} Q$ and since
$Q \in \langle J_Q \rangle$  by Euler's identity,
$L$ is in the ideal $\langle J_Q \rangle_P$ generated by the entries of
$J_Q$ in $S_P.$   We conclude that
$$
\langle J_Q \rangle_P
=
\left\langle K \frac{\partial L}{\partial x_1},
 \ldots, K \frac{\partial L}{\partial x_k}, L \right\rangle_P.
$$
Since $K$ is a unit, it can be removed.  Further, we can use Euler's relation to write the last
generator $L$ in as a linear combination of the first $k$ generators.  We thus obtain
\begin{eqnarray*}
\langle J_Q \rangle_P
& = &
\left\langle \frac{\partial L}{\partial x_1},
 \ldots, \frac{\partial L}{\partial x_k} \right\rangle_P
\\
& = &
 \langle J_{L}\rangle_P .
\end{eqnarray*}
Since localization is an exact functor (see, for example, \cite{e}),
localizing a
minimal free resolution of $\langle J_Q \rangle$ yields
a free resolution (possibly non-minimal)
of $\langle J_Q \rangle_P,$ which equals $\langle J_L \rangle_P$.
The free resolution obtained for
$\langle J_L \rangle_P$ is also a free resolution for $\langle J_L \rangle$
because $L$ is in $P.$ We conclude that
$$
\pdim(\langle J_L \rangle) \le \pdim(\langle J_Q \rangle).  
$$
This completes the proof of Theorem \ref{thm:PD}.

We remark that Theorem \ref{thm:PD} fails if one does not assume that the
subarrangement is closed.  An easy example is the braid arrangement $A_4.$
It is free but
contains three generic subarrangements of four lines, none of them closed.

Combining the results on generic arrangements and Theorem \ref{thm:PD},
we obtain a simple but useful combinatorial corollary.

\begin{cor}
\label{cor:CC}
Let $\A$ be an arrangement whose matroid
contains a generic flat of rank $r.$  Then
$\pdim(\TM) \ge r - 2.$  In particular, if the
matroid of $\A$ contains a closed circuit of size $m,$
then $\pdim(\TM) \ge m-3$.
\end{cor}

This corollary extends the folk-lore lemma that an arrangement cannot
be free if it contains a closed subarrangement consisting
of four hyperplanes in general position in $3$-dimensional space. 

Let $G$ be a graph (without loops or multiple edges) with 
vertex set $\{1,2, \ldots,n\}$ and edge set $E.$  The
{\em graphic} arrangement $\A_G$ is the collection
$\{V(x_i - x_j) \,|\, \{i,j\} \in E\}.$
For example, the braid arrangement $A_n$ is the arrangement associated to the
complete graph $K_n,$ the graph containing all possible edges.
Graphic arrangements are never essential; an arrangement from a
connected graph can be made essential by suppressing a variable.
For example, $A_n$ can be made into the essential arrangement
$\hat{A}_n$ consisting of the hyperplanes $V(x_i)$ and $V(x_i - x_j),$
where $1 \le i < j \le n-1.$

The only generic flats in graphic arrangements are closed circuits.  
Closed circuits correspond to induced cycles.
Thus, Corollary
\ref{cor:CC} also extends the reverse implication of a 
theorem (combining results in \cite{stan} and \cite{t}) that a graphic
arrangement is free if and only if its
graph is {\em chordal,} or equivalently, its graph has no induced
cycles of length greater than $3.$ In particular, we have:

\begin{cor}
\label{cor:CC1}
If a graph $G$ contains an induced cycle of length
$m$, then $\pdim(D(\A_G) \ge m-3$.
\end{cor}
 
We close this section with some illustrations of Corollary  \ref{cor:CC1}
and several related problems.

\begin{example}
\label{ex:1}
Consider the graph $G$ (with $8$ vertices) given by the $1$-skeleton of the cube:
\vskip .15in
\begin{center}
\epsfig{file=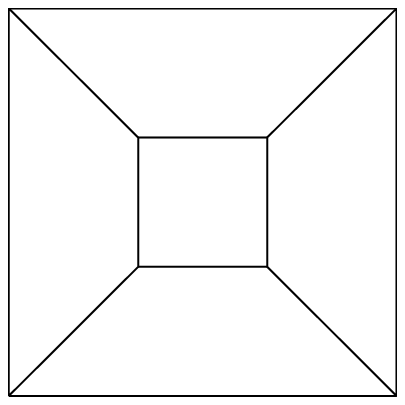,height=1.0in,width=1.2in}
\end{center}
\vskip .15in 
A free resolution for $\TMG$ is:
$$
\begin{small}
0 \longrightarrow S^3(-7) \longrightarrow S^{13}(-6) \longrightarrow 
\begin{array}{c}
 S(-4)\\
\oplus \\
S^{19}(-5)
\end{array}
\longrightarrow\\
\begin{array}{c}
S(-1)\\
 \oplus \\
S(-2)\\
 \oplus \\
S^9(-3)\\
\oplus \\
S^6(-4)
\end{array}
\longrightarrow D(A_G) \longrightarrow 0
\end{small}
$$
The diagram gives the degree (but not the explicit expressions) of the generators 
of the free modules.  For example, from the diagram, one sees that $\TMG$ can be generated by 
$17$ generators, one of degree $1$ (the Euler derivation), one of degree $2,$ nine of degree $3,$ and six of degree $4.$  
These generators have relations which can be generated by $20$ relations.  
The indexing of a free resolution starts at zero, and so $\TMG$ has projective 
dimension $3.$ Since $G$ has an induced cycle of length $6$, this
is the lower bound predicted by Corollary \ref{cor:CC1}.
\end{example}

\begin{example}
\label{ex:2}
Let $G$ be the {\em triangular prism}:
\vskip .15in
\begin{center}
\epsfig{file=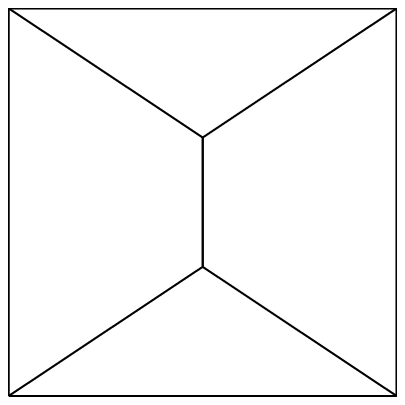,height=1.0in,width=1.2in}
\end{center}
\vskip .15in 
A free resolution for $\TMG$ is:
$$
\begin{small}
0 \longrightarrow S(-5) \longrightarrow  S^5(-4)
\longrightarrow\\
\begin{array}{c}
S(-1)\\
 \oplus \\
S(-2)\\
 \oplus \\
S^7(-3)
\end{array}
\longrightarrow D(A_G) \longrightarrow 0
\end{small}
$$
The maximum length of an induced cycle in $G$ is $4,$
but the projective dimension of $\TMG$ is $2$.  Hence,
$\pdim(\TMG)$ can
be strictly greater than the bound given in Corollary \ref{cor:CC1}.
\end{example}

Example \ref{ex:2} raises several questions.  Is there a
characterization of graphs $G$
for which $\pdim(\TMG) = m-3,$ where $m$ is the maximum size of
an induced circuit?  Are there reasonable formulas involving
graph parameters for $\pdim(\TMG)$?
In analogy to excluded minors in matroid theory (see, for example, \cite{Kung},
Section 8), define an arrangement $\A$ to be {\em $k$-minimal}
if $\pdim (D(\A)) = k$ and for every
proper closed subarrangement $\A_X \subset A,$
$\pdim (D(\A_X)) < k.$ The graphic arrangement of
the triangular prism is $2$-minimal and rank-$m$ generic arrangements
are $(m-2)$-minimal.  It seems an interesting problem to classify 
$k$-minimal arrangements.

\section {Orthogonal duals of arrangements}

Let $\A$ be a hyperplane arrangement in $n$-dimensional space.
We construct a $n \times |\A|$ matrix $A$ as follows:
each hyperplane $H$ in $\A$ labels a column equal to
$(c_1,c_2, \ldots,c_n )^t,$ where $H = V(c_1x_1 + c_2x_2 +
\ldots + c_n x_n).$
Conversely, given a matrix $A,$ we construct an arrangement by
{\em simplifying,} that is, removing
all zero columns, constructing
a multiset of hyperplanes corresponding to the kernels of the
linear forms defined by the columns, giving a multiarrangement, and 
disregarding the multiplicities to obtain an arrangement.

If $\A$ is essential, the hyperplanes in $\A$ intersect in 
the zero subspace, and the matrix $\A$ has full rank $n.$
The correspondence between essential arrangements $\A$
and $n \times |\A|$ matrices $A,$ with no zero columns and no two columns a non-zero multiple 
of each other,  
is bijective up to left multiplication by elements of $GL(n)$, 
and right multiplication by a product of a permutation matrix and a
non-singular diagonal matrix.  The matrix $A$ is a representation for
the matroid $M(\A)$ of the arrangement $\A.$

Suppose that $A$ is an $n \times d$ matrix having full rank $n.$
An {\em (orthogonal) dual} of $A$ is an $(d - n) \times d$
matrix $B$ having full rank $d - n$ such that any row of $A$
is orthogonal (under the usual dot product) to any row
of $B.$  The
matrix $B$ exists and is determined
to left multiplication
by a non-singular matrix.
In addition, if $A^\perp$ is a dual of $A,$ then it is also a
dual of any matrix obtained from
$A$ by left multiplication by a non-singular matrix.  Thus, duality is
an operation defined between equivalence classes of matrices.  In particular,
there is an easy way to construct a dual of $A.$   Put $A$ into the form
$[I|C],$ where $I$ is the $n \times n$ identity matrix.  Then
a dual of $A$ is $[-C^t|I],$ where $I$ is the $(d-n) \times (d-n)$ identity matrix.

If $\A$ is an essential
arrangement with matrix $A,$ we define its
{\em (orthogonal) dual} $\A^\perp$ to be the arrangement obtained
from a dual of the matrix $A.$  Note that
because we discard zero columns and ignore multiplicities,
$\A$ is not reconstructible from $\A^\perp$ in general.

We will also need several elementary facts from the theory of matroid duality
(see \cite{CR, Crapo, OX, Abs}).
There are many ways to define the orthogonal dual of a matroid.  For us, the best
definition is the circuit-cocircuit definition.  Recall that
a {\em circuit} is a minimal dependent set and a
{\em cocircuit} is the complement of a copoint.
The {\em (orthogonal) dual} $M^\perp$ of $M$
is the matroid on the same ground set
whose circuits are exactly the cocircuits of $M.$
Duality interchanges contraction and deletion,
that is, for a subset $B$ of the set of elements, $(M^{\perp})/B$ equals
$(M \backslash B)^\perp.$
It is true (but not obvious) that the matroid of the dual
arrangement $\A^\perp$
is the simplification of the dual of the matroid of $\A.$
Despite its age, the
neatest and most accessible proof of this
remains Whitney's original proof in \cite{Abs}.

A {\em loop} is an element $e$ such that
the set $\{e\}$ is a circuit.  An {\em isthmus} is an element
$e$ such that $\{e\}$ is a cocircuit, so that $M$ is the direct sum
$(M \backslash \{e\}) \oplus \{e\}.$
For graphs, an isthmus is an edge whose removal increases the number
of connected components.
Duality interchanges loops and isthmuses.

We shall call closure in the dual matroid $M^{\perp}$
{\em $\perp$-closure.}

\begin{lem}
\label{lem:Perp-closure}
Let $M$ be a matroid on the set $E$ and $B \subseteq E.$
Then $e$ is in the $\perp$-closure of $B$ if and only if
$e$ is in $B$ or $e$ is an isthmus in the deletion $M \backslash B.$
In particular, a cocircuit $B$ is
$\perp$-closed if its complementary copoint $X$ has no isthmuses.
\end{lem}

\begin{proof}
The lemma follows from dualizing the statement:
a point $e$ is in the $\perp$-closure
if and only if $e$ is a loop
in the contraction $M^\perp / B.$
\end{proof}

The {\em cycle matroid} $M(G)$ of a graph is the matroid on the edge set
whose circuits are the cycles of the graph.
A {\em cutset} in a graph $G$ is an edge-subset whose removal increases
the number of connected components of $G.$
The circuits of the dual matroid $G^{\perp}$ are precisely the
minimal cutsets of $G.$
For graphs, an isthmus is an edge which is a cutset by itself.
Thus, Lemma \ref{lem:Perp-closure} gives an
easy way to determine whether a minimal cutset is $\perp$-closed.
We remark that the
set of all edges incident on a vertex $v$ is a minimal cutset.
Such ``vertex cutsets''
usually contain few edges compared
to other minimal cutsets.

The complete graph $K_n$ is the graph on $n$ vertices with all
possible edges.
The maximum size of a minimal cutset in complete graphs is given
in the next lemma.

\begin{thm}
\label{thm:KN}
 The largest cocircuit
in the cycle matroid $M(K_{n+1})$ has size $\lceil n(n+2)/4 \rceil.$
When $n \ge 5,$ the largest cocircuits are $\perp$-closed.
\end{thm}
\begin{proof}
Because $K_{n+1}$ contains all possible edges, minimal cutsets are
in bijection with partitions of the vertex set into
two non-empty subsets and these cutsets disconnect $K_{n+1}$ into
two disjoint smaller complete graphs.
By Lemma \ref{lem:Perp-closure}, every minimal cutset
in $K_{n+1}$ gives a $\perp$-closed cocircuit
with the exception of the cutsets which divide
$K_{n+1}$ into a $K_{n-1}$ and a single edge $K_2.$
The largest minimal cutsets
are those which divide $K_{n+1}$ into two connected components of almost
equal size.  We conclude that the largest cocircuit
in $M(K_{n+1})$ has size $k^2$ if $n + 1 = 2k$ and
$k(k+1)$ if $n+1 = 2k + 1.$
To finish, it is easy to check that
$\lceil n(n+2)/4 \rceil$ equals $k^2$ or $k(k+1)$
depending on the parity of $n+1.$
\end{proof}

\begin{thm}
\label{thm:Mod}
Let $M$ be a rank-$n$ matroid on the set $S$ with a modular copoint $X.$
Suppose that the cocircuit $S \backslash X$ spans.   If $n \ge 5,$
then there exists a $\perp$-closed cocircuit in $M$ of size
at least $\lceil n(n+2)/4 \rceil.$
\end{thm}
\begin{proof}
We shall use the following lemma.

\begin{lem}
\label{lem:Cat}
Under the hypotheses in the theorem, $M$
 contains a spanning submatroid isomorphic to
$M(K_{n+1}).$
\end{lem}
\begin{proof}
This is a combination of Lemma 5.3 in \cite{Num} and Lemma 5.14 in \cite{Kung}.
For the sake of completeness, we will give a proof in the language of
arrangements and linear forms.  Choose coordinates so that
the linear forms $x_i, 1 \le i \le n$ are in the cocircuit $S \backslash X$
and the copoint $X$ is the subarrangement of all linear forms whose
kernel contains the point $(1,1, \ldots,1).$  By modularity,
$$
\rk((x_i \vee x_j) \wedge X) = \rk(X) + \rk(x_i \vee x_j) - n = 1
$$
for every pair $x_i$
and  $x_j$ of linear forms.  Hence,
there is a linear combination of $x_i$ and $x_j$ whose kernel
contains $(1,1,\ldots,1).$  This form is $x_i - x_j,$  so 
the arrangement contains the subarrangement
$\{x_i, \, x_i - x_j \, | \, 1 \le i < j \le n\},$ which is
 isomorphic to the graphic arrangement of $K_{n+1}.$
\end{proof}

Let $K$ be a spanning submatroid in $M$ isomorphic to $M(K_{n+1}).$
Take a copoint $X^\prime$ in the submatroid $K$ and let $X$ be the
closure of $X^{\prime}$ in $M.$  A point in
$K \backslash X^{\prime}$ is still not in $X.$   Hence,
$$
|S \backslash X| \ge |K \backslash X^{\prime}|.
$$
Choosing $X^\prime$ in $M|K$ so that $K \backslash X^{\prime}$
has size $\lceil n(n+2)/4 \rceil,$
we obtain a cocircuit in $M$ having size at least $\lceil n(n+2)/4 \rceil.$
If $n \ge 5,$ the copoint $X^\prime$ in $K$ contains no isthmuses.
Since $X$ and $X^\prime$ have the same rank,
a direct summand of $X$ induces a direct summand of $X^{\prime}.$
As $X^\prime$ contains no isthmuses,
$X$ also contains no isthmuses and the cocircuit $S \backslash X$
is $\perp$-closed.
\end{proof}

Corollary \ref{cor:CC} and Theorem \ref{thm:Mod} imply
Theorem \ref{thm:dualPD}.

Since the matroid of the braid arrangement $A_n$ is
$M(K_n),$ Theorem \ref{thm:PD} and Theorem \ref{thm:KN} imply
that the projective dimension of the dual of the ``essential'' braid
arrangement $\hat{A}_n$ is at least $\lceil (n-1)(n+1)/4 \rceil - 3.$
Lower bounds for the other families of real reflection arrangements
can be obtained using the method in the proof of Theorem \ref{thm:KN}.

\begin{thm}
\label{thm:reflection}
When $n \ge 5,$
$$
\pdim (B_n^\perp) \ge \lfloor \tfrac {2}{3} n^2 +
\tfrac {1}{3}n - \tfrac {1}{24} \rfloor - 3.
$$
When $n \ge 6,$
$$
\pdim (D_n^\perp) \ge \lfloor \tfrac {2}{3} n^2 -
\tfrac {1}{3}n + \tfrac {1}{24} \rfloor - 3.
$$
\end{thm}
\begin{proof}
Consider
the copoint isomorphic to the direct sum $A_k \oplus B_{n-k}$
in $B_n$ spanned by the $n-1$ linear forms
$$
x_1 - x_2, x_2 - x_3, \ldots, x_{k-1} - x_k, x_{k+1}, x_{k+2}, \ldots, x_n.
$$
The cocircuit complementary to $X$ has size
$$
{k \choose 2} +  k + 2k(n-k).
\eqno(1)$$
We obtain a cocircuit
of maximum size when $k$ is the integer closest to $2n/3 + \tfrac {1}{6}$
and this maximum size is obtained by substitution into formula (1) and
rounding down.
Since $A_2$ and $B_1$ contain a single form,
the cocircuits of maximum size are $\perp$-closed if $n \ge 5.$
The argument for $D_n$ is similar.
\end{proof}

The argument for $B_n$ can also be applied to the complex
reflection arrangements $G(n,1,l)$ to give a rough lower bound of
$(l^2/(2l+2)) n^2$ for $\pdim (G(n,1,l)^\perp)$
when $n \ge 5.$

\vskip 0.3in
\noindent{\bf Acknowledgment} 
The Macaulay2 software package, available at 
\[
{\tt http://www.math.uiuc.edu/Macaulay2/}
\]
allowed us to compute many examples which provided evidence for the
paper.  In particular, Examples \ref{ex:1} and \ref{ex:2} were computed using this
software.

\bibliographystyle{amsalpha}

\end{document}